\documentstyle[amsfonts,12pt]{article}

\normalbaselines
\input{tcilatex}

\begin{document}

\title{Restricted linear properties of scalar product in T-geometry}
\author{Yuri A.Rylov}
\date{Institute for Problems in Mechanics, Russian Academy of Sciences,\\
101, bild.1 Vernadskii Ave., Moscow, 119526, Russia. \\
e-mail: rylov@ipmnet.ru}
\maketitle

\begin{abstract}
It is shown that scalar product of two vectors can be introduced in any 
geometry (metric space) independently of possibility of the linear space 
introduction. In general, linear properties of scalar product are restricted. 
Domain of linearity increases in homogeneous geometry.
\end{abstract}

T-geometry \cite{R02} is a very general geometry giving on an arbitrary set $%
\Omega $ of points by setting the world function $\sigma $: 
\begin{equation}
\sigma :\;\;\;\Omega \times \Omega \rightarrow {\Bbb R},\qquad \sigma \left(
P,P\right) =0,\qquad \forall P,Q\in \Omega  \label{a1}
\end{equation}
T-geometry is symmetric, if additionally 
\begin{equation}
\sigma \left( P,Q\right) =\sigma \left( Q,P\right)  \label{a1.1}
\end{equation}
The T-geometry has no its own axiomatics. It is obtained as a result of a
deformation of the proper Euclidean geometry. In other words, any relation
of T-geometry is obtained as a result of some deformation of the
corresponding relation of the proper Euclidean geometry. The relation of
Euclidean geometry is written in $\sigma $-immanent form (i.e. in terms of
only world function $\sigma $). Thereafter the world function of the proper
Euclidean space is replaced by the world function of the T-geometry in
question. Corresponding relation of T-geometry arises. Representation of a
proper Euclidean relation in $\sigma $-immanent form is always possible,
because the proper Euclidean geometry can be formulated $\sigma $-immanently
(in terms of only world function) \cite{R02}.

Vector ${\bf P}_{0}{\bf P}_{1}\equiv \overrightarrow{P_{0}P_{1}}$ in
T-geometry is the ordered set of two points ${\bf P}_{0}{\bf P}_{1}=\left\{
P_{0},P_{1}\right\} $, \ \ $P_{0},P_{1}\in \Omega $. The scalar product $%
\left( {\bf P}_{0}{\bf P}_{1}.{\bf Q}_{0}{\bf Q}_{1}\right) $ of two vectors 
${\bf P}_{0}{\bf P}_{1}$ and ${\bf Q}_{0}{\bf Q}_{1}$ is defined by the
relation 
\begin{equation}
\left( {\bf P}_{0}{\bf P}_{1}.{\bf Q}_{0}{\bf Q}_{1}\right) =\sigma \left(
P_{0},Q_{1}\right) +\sigma \left( P_{1},Q_{0}\right) -\sigma \left(
P_{0},Q_{0}\right) -\sigma \left( P_{1},Q_{1}\right) ,  \label{a2}
\end{equation}
for all $P_{0},P_{1},Q_{0},Q_{1}\in \Omega $. As it follows from (\ref{a1})
-- (\ref{a2}), in the symmetric T-geometry 
\begin{equation}
\left( {\bf P}_{0}{\bf P}_{1}.{\bf Q}_{0}{\bf Q}_{1}\right) =\left( {\bf Q}%
_{0}{\bf Q}_{1}.{\bf P}_{0}{\bf P}_{1}\right) ,\qquad
P_{0},P_{1},Q_{0},Q_{1}\in \Omega  \label{a2.1}
\end{equation}

When the world function $\sigma $ is such one \cite{R02} that the $\sigma $%
-space $V=\left\{ \sigma ,\Omega \right\} $ is the $n$-dimensional proper
Euclidean space $E_{n}$ the scalar product (\ref{a2}) turns to scalar
product of two vectors in $E_{n}$. Besides as it follows from (\ref{a1}), (%
\ref{a2}) that in any T-geometry 
\begin{equation}
\left( {\bf P}_{0}{\bf P}_{1}.{\bf Q}_{0}{\bf Q}_{1}\right) =-\left( {\bf P}%
_{0}{\bf P}_{1}.{\bf Q}_{1}{\bf Q}_{0}\right) ,\qquad \forall
P_{0},P_{1},Q_{0},Q_{1}\in \Omega  \label{a3}
\end{equation}
\begin{equation}
\left( {\bf Q}_{0}{\bf Q}_{1}.{\bf P}_{0}{\bf P}_{1}\right) =-\left( {\bf Q}%
_{0}{\bf Q}_{1}.{\bf P}_{1}{\bf P}_{0}\right) ,\qquad \forall
P_{0},P_{1},Q_{0},Q_{1}\in \Omega  \label{a4}
\end{equation}
\begin{eqnarray}
\left( {\bf P}_{0}{\bf P}_{1}.{\bf Q}_{0}{\bf Q}_{1}\right) +\left( {\bf P}%
_{1}{\bf P}_{2}.{\bf Q}_{0}{\bf Q}_{1}\right) &=&\left( {\bf P}_{0}{\bf P}%
_{2}.{\bf Q}_{0}{\bf Q}_{1}\right) ,  \label{a5} \\
\left( {\bf Q}_{0}{\bf Q}_{1}.{\bf P}_{0}{\bf P}_{1}\right) +\left( {\bf Q}%
_{0}{\bf Q}_{1}.{\bf P}_{1}{\bf P}_{2}\right) &=&\left( {\bf Q}_{0}{\bf Q}%
_{1}.{\bf P}_{0}{\bf P}_{2}.\right) ,  \label{a5.0}
\end{eqnarray}
for all $P_{0},P_{1},P_{2},Q_{0},Q_{1}\in \Omega $

Let us consider relations between two vectors ${\bf P}_{0}{\bf P}_{1}$ and $%
{\bf R}_{0}{\bf R}_{1}$%
\begin{equation}
{\bf P}_{0}{\bf P}_{1}{\cal T}{\bf R}_{0}{\bf R}_{1}:\quad \left( {\bf P}_{0}%
{\bf P}_{1}.{\bf Q}_{0}{\bf Q}_{1}\right) =\left( {\bf R}_{0}{\bf R}_{1}.%
{\bf Q}_{0}{\bf Q}_{1}\right) \wedge \left( {\bf Q}_{0}{\bf Q}_{1}.{\bf P}%
_{0}{\bf P}_{1}\right) =\left( {\bf Q}_{0}{\bf Q}_{1}.{\bf R}_{0}{\bf R}%
_{1}\right)  \label{a5.1}
\end{equation}
for all $Q_{0},Q_{1}\in \Omega $. The relation (\ref{a5.1}) is reflexive,
symmetric and transitive. It can be considered to be equivalence relation.
Thus, two vectors ${\bf P}_{0}{\bf P}_{1}$ and ${\bf R}_{0}{\bf R}_{1}$ are
equivalent $\left( {\bf P}_{0}{\bf P}_{1}={\bf R}_{0}{\bf R}_{1}\right) $,
if the conditions 
\begin{eqnarray}
\left( {\bf P}_{0}{\bf P}_{1}.{\bf Q}_{0}{\bf Q}_{1}\right) &=&\left( {\bf R}%
_{0}{\bf R}_{1}.{\bf Q}_{0}{\bf Q}_{1}\right) ,\qquad \forall Q_{0},Q_{1}\in
\Omega  \label{a5.2} \\
\left( {\bf Q}_{0}{\bf Q}_{1}.{\bf P}_{0}{\bf P}_{1}\right) &=&\left( {\bf Q}%
_{0}{\bf Q}_{1}.{\bf R}_{0}{\bf R}_{1}\right) ,\qquad \forall Q_{0},Q_{1}\in
\Omega  \label{a5.3}
\end{eqnarray}
take place. Then one obtains from (\ref{a3}) -- (\ref{a5.0}) 
\begin{equation}
{\bf P}_{1}{\bf P}_{0}=-{\bf P}_{0}{\bf P}_{1}=\alpha {\bf P}_{0}{\bf P}%
_{1},\qquad \alpha =-1  \label{a6}
\end{equation}
\begin{equation}
{\bf P}_{0}{\bf P}_{2}={\bf P}_{0}{\bf P}_{1}+{\bf P}_{1}{\bf P}_{2},\qquad 
{\bf P}_{0}{\bf P}_{1}={\bf P}_{0}{\bf P}_{2}-{\bf P}_{1}{\bf P}_{2}
\label{a7}
\end{equation}

Relations (\ref{a6}), (\ref{a7}) determine multiplication of a vector by a
real number $\alpha =\pm 1$ and summation of vectors. But vectors in the $%
\sigma $-space $V=\left\{ \sigma ,\Omega \right\} $ are not free vectors and
two vectors may be added, provided end of one vector is the origin of other
vector. At such a definition of the vector summation the scalar product of
vectors have linear properties 
\begin{eqnarray}
\left( \alpha {\bf P}_{0}{\bf P}_{1}{\bf +\beta R_{0}R_{1}.Q}_{0}{\bf Q}%
_{1}\right) &=&\alpha \left( {\bf P}_{0}{\bf P}_{1}{\bf .Q}_{0}{\bf Q}%
_{1}\right) +{\bf \beta }\left( {\bf R_{0}R_{1}.Q}_{0}{\bf Q}_{1}\right)
\label{a8} \\
\left( {\bf Q}_{0}{\bf Q}_{1}.\alpha {\bf P}_{0}{\bf P}_{1}{\bf +\beta
R_{0}R_{1}}\right) &=&\alpha \left( {\bf Q}_{0}{\bf Q}_{1}.{\bf P}_{0}{\bf P}%
_{1}\right) +{\bf \beta }\left( {\bf Q}_{0}{\bf Q}_{1}.{\bf R_{0}R_{1}}%
\right)  \label{a8.1}
\end{eqnarray}
\[
\alpha =\pm 1,\beta =\pm 1,\qquad P_{0},P_{1},R_{0},R_{1}\in \Omega
,\;\;\forall Q_{0},Q_{1}\in \Omega 
\]

The linear properties of the scalar product are restricted as compared with
those in the Euclidean space, because according to (\ref{a3}) -- (\ref{a5.0}%
) the linear operations on vectors are always defined only in the following
cases:

\begin{enumerate}
\item  $\alpha =0\vee \beta =0$

\item  $\alpha =\beta =\pm 1,$ $\left( P_{1}=R_{0}\right) \vee \left(
R_{1}=P_{0}\right) $

\item  $\alpha =-\beta =\pm 1,$ $\left( P_{1}=R_{1}\right) \vee \left(
P_{0}=R_{0}\right) $
\end{enumerate}

In other cases the linear operations on vectors are not defined, in general,
but they may be defined in some special T-geometries, when relations (\ref
{a8}), (\ref{a8.1}) 
\begin{equation}
\begin{array}{c}
\left( {\bf S_{0}S_{1}.Q}_{0}{\bf Q}_{1}\right) =\alpha \left( {\bf P}_{0}%
{\bf P}_{1}{\bf .Q}_{0}{\bf Q}_{1}\right) +{\bf \beta }\left( {\bf %
R_{0}R_{1}.Q}_{0}{\bf Q}_{1}\right) \\ 
\left( {\bf Q}_{0}{\bf Q}_{1}.{\bf S_{0}S_{1}}\right) =\alpha \left( {\bf Q}%
_{0}{\bf Q}_{1}.{\bf P}_{0}{\bf P}_{1}\right) +{\bf \beta }\left( {\bf Q}_{0}%
{\bf Q}_{1}.{\bf R_{0}R_{1}}\right)
\end{array}
\label{a9.2}
\end{equation}
take place for all vectors ${\bf Q_{0}Q_{1}}$ and some $\alpha ,\beta \in 
{\Bbb R}$. Then by definition 
\begin{equation}
{\bf S_{0}S_{1}}=\alpha {\bf P}_{0}{\bf P}_{1}+\beta {\bf R_{0}R_{1}}
\label{a10}
\end{equation}
In general case, there is no vector ${\bf S_{0}S_{1}}$, satisfying (\ref
{a9.2}), for all vectors ${\bf Q_{0}Q_{1}}$, and linear operations are
defined only in the case, when (\ref{a3}) -- (\ref{a5.0}) are fulfilled.

Let $E_{n}=\left\{ \sigma _{{\rm E}},{\Bbb R}^{n}\right\} $ be $n$%
-dimensional proper Euclidean space. Remove some set $\Omega ^{\prime }$ of
points from ${\Bbb R}^{n}$. The $\sigma $-space $V=\left\{ \sigma _{{\rm E}%
},\Omega \right\} $, \ $\Omega ={\Bbb R}^{n}\backslash \Omega ^{\prime }$ is
a result of deformation of Euclidean space $E_{n}$, because deletion of
points is a special kind of deformation. Linear operation (\ref{a10}) is
defined in $E_{n}$ for all points $P_{0},P_{1},R_{0},R_{1}$ and all real $%
\alpha ,\beta $. Linear operation (\ref{a10}) is defined in $V$ for some
real $\alpha ,\beta $ and for $P_{0},P_{1},R_{0},R_{1}\in \Omega $, provided
there are such points $S_{0},S_{1}\in \Omega $ that relations (\ref{a9.2})
are fulfilled for all $Q_{0},Q_{1}\in \Omega $. Thus, deletion of points
from Euclidean space destroys its linear structure, but not completely. A
change of the world function plays the same role. Nevertheless, there is a
minimal linear structure (\ref{a6}), (\ref{a7}) which remains at any
deformation of the Euclidean space.

\end{document}